\documentclass{article}
\usepackage{amsmath,amsthm,amssymb}

\newtheorem{clm}{Claim}
\newtheorem{thm}{Theorem}
\newtheorem{lem}{Lemma}
\newtheorem{coro}{Corollary}
\theoremstyle{definition}
\newtheorem{defn}{Definition}
\newtheorem{prop}{Proposition}

\begin{document}
\title{Local Ramsey theory. An abstract approach}
\author{Mijares Jos\'e\thanks{jmijares@euler.ciens.ucv.ve}\\Escuela de Matem\'atica\\Universidad Central de Venezuela \and Nieto Jes\'us\thanks{jnieto@usb.ve}\\Depto. Formaci\'on General y Ciencias B\'asicas
\\Universidad Sim\'on Bol\'ivar}

\maketitle

\begin{abstract}
It is shown that the known notion of selective coideal can be extended to a family $\mathcal{H}$ of subsets of $\mathcal{R}$, where $(\mathcal{R},\leq,r)$ is a topological Ramsey space in the sense of Todorcevic (see \cite{todo}). Then it is proven that, if $\mathcal{H}$ selective, the $\mathcal{H}$-Ramsey and $\mathcal{H}$-Baire subsets of $\mathcal{R}$ are equivalent. This extends the results of Farah in \cite{farah} for semiselective coideals of $\mathbb{N}$. Also, it is proven that the family of ${\cal H}$--Ramsey subsets of ${\cal R}$ is closed under the Souslin operation.
\end{abstract}

\section{Introduction}
In \cite{mathias}, Mathias introduces the \emph{happy families} (or selective coideals) of subsets of $\mathbb{N}$ and relativizes the notion of completely Ramsey (see \cite{galpri}) subsets of $\mathcal{P}(\mathbb{N})$ (the set of subsets of $\mathbb{N}$) to such families. Then he proves that analitic sets are $\mathcal{F}$-Ramsey when $\mathcal{F}$ is a Ramsey ultrafilter and generalizes this result for arbitrary happy families. In \cite{farah}, Farah gives an answer to the question of Todorcevic: what are the combinatorial properties of the family $\mathcal{H}$ of ground model subsets of $\mathbb{N}$ which warranties diagonalization of the Borel partitions? This is done by imposing a condition on $\mathcal{H}$ which is weaker than selectivity, that is the notion of \emph{semiselectivity}. In that work he proves that the semiselectivity of $\mathcal{H}$ is enough for a subset of $\mathbb{N}^{[\infty]}$ to be $\mathcal{H}$-Ramsey if and only if it has the (abstract) Baire property with respect to $\mathcal{H}$. In \cite{mij2}, Mijares extends this result to any \emph{topological Ramsey space} (see \cite{todo}) by generalizing the notion of Ramsey ultrafilter to such spaces. In this work, it is proven that a family $\mathcal{H}$ of subsets of a topological Ramsey space $\mathcal{R}$, provided with suitable features, corresponds to the semiselective coideal given by Farah in the above mentioned work. Then the results about $\mathcal{H}$--Ramsey and $\mathcal{H}$--Baire sets of $\mathcal{R}$ are extended to this context.

\vspace{0.25cm}

\noindent The structure of this work is as follows: In section \ref{topo}, some material regarding the so called \emph{topological Ramsey theory} is given, see \cite{todo} or \cite{carsimp}. In section \ref{theory} it is proven that certain features of a family of sets of $\mathcal{R}$ are sufficient for a subset of $\mathcal{R}$ to be $\mathcal{H}$--Ramsey if ad only if it is $\mathcal{H}$--Baire. This is done by defining the \emph{D--O property} which is proposed as the corresponding notion to  dense open sets. This is what is known as \emph{local Ramsey theory}. In section 4, it is shown that the family of ${\cal H}$--Ramsey subsets of ${\cal R}$ is closed under the Souslin operation if ${\cal H}$ is selective. Finally, in section \ref{examples}, examples for which the results hold are given. 

\section{Preliminaries: Topological Ramsey Theory.}\label{topo}

The definitions and results throughout this section  are
expected to appear in \cite{todo}. A previous presentation of the
following notions can also be found in \cite{carsimp}.

\vspace{0.25cm}

\noindent Consider a triplet of the form $(\mathcal{R}, \leq,
r )$, where $\mathcal{R}$ is a set, $\leq$
is a quasi order on $\mathcal{R}$, $\mathbb{N}$ is
 the set of natural numbers, and $r: \mathbb{N}\times\mathcal{R}\rightarrow \mathcal{AR}$
 is a function with range $\mathcal{AR}$. For each $A\in \mathcal{R}$, we say that $r_{n}(A)$ is
 \textit{the} $n$th \textit{approximation of} $A$. Denote, for every $n\in \mathbb{N}$ and every $A\in \mathcal{R}$, $r(n,A)=r_n(A)$ and $r_n(\mathcal{R})=\mathcal{AR}_n$. In order to capture the combinatorial structure
 required to ensure the provability of an Ellentuck type theorem, some assumptions on
 $(\mathcal{R}, \leq, r)$ will be imposed. The first three of them are the following:

\begin{itemize}
\item[{(A.1)}]For any $A\in \mathcal{R}$, $r_{0}(A) = \emptyset$.
\item[{(A.2)}]For any $A,B\in \mathcal{R}$, if $A\neq B$ then
$(\exists n)r_{n}(A)\neq r_{n}(B)$. 
\item[{(A.3)}]If $r_{n}(A) =
r_{m}(B)$ then $n = m$ and $(\forall i<n)r_{i}(A) = r_{i}(B)$.
\end{itemize}
These three assumptions allow us to identify each $A\in
\mathcal{R}$ with the sequence $(r_{n}(A))_{n}$ of its
approximations. In this way, if $\mathcal{AR}$ has the discrete
topology, $\mathcal{R}$ can identified with a subspace of the
(metric) space $\mathcal{AR}^{[\infty]}$ (with the product
topology) of all the sequences of elements of $\mathcal{AR}$. Via
this identification, $\mathcal{R}$ will be regarded as a subspace of
$\mathcal{AR}^{[\infty]}$, and we will say that $\mathcal{R}$ is
{\it metrically closed} if it is a closed subspace of
$\mathcal{AR}^{[\infty]}$.

\vspace{.25 cm}

\noindent Also, for $a\in \mathcal{AR}$, define the \textit{length} of $a$,
$|a|$, as the unique $n$ such that $a = r_{n}(A)$ for some $A\in
\mathcal{R}$, and the \textit{Ellentuck type
neighborhoods} on $\mathcal{R}$ 
$$[a,A] = \{B\in \mathcal{R} : (\exists n)(a =
r_{n}(B))\ \ \mbox{and} \ \ B\leq A\}$$ where $a\in\mathcal{AR}$
and $A\in \mathcal{R}$. If $[a,A]\neq\emptyset$ we will say that
$a$ is \textit{compatible} with $A$ (or $A$ is compatible with
$a$). Let $\mathcal{AR}[A] = \{a\in \mathcal{AR} : a\ \ \mbox{is
compatible with}\ \ A\}$.

\vspace{.25 cm}

\noindent Denote $[n,A]$ for $[r_{n}(A),A]$, and $Exp(\mathcal{R})$ for
the family of all the neighborhoods $[n,A]$. This family generates
the natural "exponential" topology on $\mathcal{R}$ which is finer
than the product topology.

\vspace{.25 cm}

\noindent Now, an analog notion for subsets of
$\mathcal{R}$, to that of \textit{Ramseyness} for subsets of
$\mathbb{N}^{[\infty]}$ is defined:

\begin{defn}
A set $\mathcal{X}\subseteq \mathcal{R}$
is \textbf{Ramsey} if for every neighborhood $[a,A]\neq\emptyset$
there exists  $B\in [a,A]$ such that $[a,B]\subseteq \mathcal{X}$
or $[a,B]\cap \mathcal{X} = \emptyset$. A set
$\mathcal{X}\subseteq \mathcal{R}$ is \textbf{Ramsey null} if for
every neighborhood $[a,A]$ there exists  $B\in [a,A]$ such that
$[a,B]\cap \mathcal{X} = \emptyset$.
\end{defn}

\begin{defn}
We say that $(\mathcal{R}, \leq,
r)$ is a \textbf{Ramsey space} if subsets of
$\mathcal{R}$ with the Baire property are Ramsey and meager
subsets of $\mathcal{R}$ are Ramsey null.
\end{defn}

\vspace{.25 cm}

\noindent In \cite{todo} it is shown that (A.1), (A.2) and (A.3), together with the
following three assumptions are conditions of sufficiency for a
triplet $(\mathcal{R}, \leq, r)$, with $\mathcal{R}$
metrically closed, to be a Ramsey space:

\vspace{.25 cm}

\noindent(A.4)(\textit{Finitization}) There is a quasi order $\leq_{fin}$ on
$\mathcal{AR}$ such that:
    \begin{itemize}
    \item[{(i)}]$A\leq B$ iff
    $\forall n\exists m \ \ r_{n}(A)\leq_{fin} r_{m}(B)$.
    \item[{(ii)}]$\{b\in \mathcal{AR} : b\leq_{fin} a\}$ is finite, for every
    $a\in \mathcal{AR}$.
    \end{itemize}

\vspace{.25 cm}

\noindent Given $a$ and $A$, we define the \textit{depth of} $a$ \textit{in}
$A$, $depth_{A}(a)$, as the minimal $n$ such that $a\leq_{fin}
r_{n}(A)$.

\begin{itemize}

\item[{(A.5)}](\textit{Amalgamation}) Given compatible $a$ and $A$
with $depth_{A}(a)=n$, the fo\-llowing holds:
\begin{itemize}
\item[{(i)}] $\forall B\in [n,A]\ \ ([a,B]\neq\emptyset)$.
\item[{(ii)}] $\forall B\in [a,A]\ \ \exists A'\in [n,A]\ \
([a,A']\subseteq [a,B])$.
\end{itemize}
\item[{(A.6)}](\textit{Pigeon Hole Principle}) Given compatible $a$
and $A$ with $depth_{A}(a) = n$, for each partition
$\phi:\mathcal{AR}_{|a|+1}\rightarrow \{0,1\}$ there is $B\in
[n,A]$ such that $\phi$ is constant in $r_{|a|+1}[a,B]$.
\end{itemize}

\vspace{0.4cm}

\noindent{\bf Abstract Ellentuck Theorem:}

\begin{thm}[Carlson]\label{AbsEll}
Any $(\mathcal{R}, \leq, r)$
with $\mathcal{R}$ metrically closed and satisfying (A.1)-(A.6) is a
Ramsey space. 
\end{thm} 

\qed

\vspace{.25 cm}

\noindent For instance, take  $\mathcal{R} = \mathbb{N}^{[\infty]}$, the set
of infinite subsets of $\mathbb{N}$, $\leq \ \ = \ \ \subseteq$
and $r_{n}(A)$ = the first $n$ elements of $A$, for each $A\in
\mathbb{N}^{[\infty]}$ . So, the set of approximations is
$\mathcal{AR} = \mathbb{N}^{[<\infty]}$, the set of finite subsets
of $\mathbb{N}$. The family of neighborhoods $[a,A]$, with
$a\in\mathbb{N}^{[<\infty]}$ and $A\in \mathbb{N}^{[\infty]}$, is
the family of Ellentuck neighborghoods. Define $\leq_{fin}$ as $a\leq_{fin} b$ iff
($a=b=\emptyset$ or $a\subseteq b$ and $max(a)=max(b)$), for
$a,b\in\mathbb{N}^{[<\infty]}$. With these definitions, (A.1)-(A.6) hold. In this case (A.6) reduces to a natural variation of
the classical pigeon hole principle for finite partitions of an
infinite set of natural numbers. Note also that
$\mathbb{N}^{[\infty]}$ is easily identified with a closed
subspace of $\mathcal{AR}^{[\infty]}$, namely, the set of all the
sequences $(x_{n})_{n}$ of finite sets such that $x_{n} =
x_{n+1}\setminus \{max(x_{n+1})\}$, for each $n\in\mathbb{N}$.
Then $(\mathbb{N}^{[\infty]}, \subseteq, r)$ is a
Ramsey space in virtue of the abstract Ellentuck theorem. Hence,
Ellentuck's theorem is obtained as a corollary:

\vspace{.25 cm}

\begin{coro}[Ellentuck]\label{Ell}
Given $\mathcal{X}\subseteq\mathbb{N}^{[\infty]}$, the
following hold:
\begin{itemize}
\item[{(a)}]$\mathcal{X}$ is Ramsey iff $\mathcal{X}$ has the
Baire Property, relative to Ellentuck's topology.
\item[{(b)}]$\mathcal{X}$ is Ramsey null iff $\mathcal{X}$ is
meager, relative to Ellentuck's topology. 
\end{itemize}
\end{coro}

\qed

\section{Selectivity}\label{theory}

\noindent From now on suppose that $({\cal R}, \leq, r)$ is a topological Ramsey space; that is, (A.1)$\dots$(A.6) hold and ${\cal R}$ is metrically closed. The following features are inspired on the known notion of coideal, so it will be used the same name: we say that ${\cal H}\subseteq {\cal R}$ is a \textbf{\emph{coideal}} of $\mathcal{R}$ (or simply a \emph{coideal}) if it satisfies:
\begin{enumerate}
\item If $A\leq B$ and $A\in {\cal H}$ then $B\in {\cal H}$.
\item \textbf{(A.5) mod ${\cal H}$:} Given $A\in {\cal H}$ and $a\in {\cal AR}(A)$, if $depth_A(a)=n$, then:
\begin{itemize}
\item[i)] $\forall B\in [n,A]\cap {\cal H}$ $([a,B]\cap {\cal H}\neq \emptyset)$.
\item[ii)] $\forall B\in [a,A]\cap {\cal H}$ $\exists A'\in [n,A]\cap {\cal H}$ $([a,A']\subseteq [a,B])$.
\end{itemize}
\item \textbf{(A.6) mod ${\cal H}$:} Given $a\in {\cal AR}$ with length $l$ and ${\cal O}\subseteq {\cal AR}_{l+1}$. Then, for every $A\in {\cal H}$ with $[a,A]\neq \phi$, there exists $B\in[depth_A(a), A]\cap {\cal H}$ such that $r_{l+1}([a,B])\subseteq {\cal O}$ or $r_{l+1}([a,X])\subseteq {\cal O}^{\,c}$.
\end{enumerate}

\noindent The natural definitions of ${\cal H}$-Ramsey and ${\cal H}$-Baire sets will be:

\begin{defn} ${\cal X}\subseteq {\cal R}$ is \textbf{${\cal H}$-Ramsey} if for every $[a,A]\neq \emptyset$, with $A\in {\cal H}$, there exists $B\in [a,A]\cap {\cal H}$ with $[a,B]\neq \emptyset$ such that $[a,B]\subseteq {\cal X}$ or $[a,B]\subseteq {\cal X}^c$. If for every $[a,A]\neq \emptyset$, there exists $B\in [a,A]\cap{\cal H}$ with $[a,B]\neq \emptyset$ such that $[a,B]\subseteq {\cal X}^c$; we say that  ${\cal X}$ is \textbf{${\cal H}$-Ramsey null}.
\end{defn}

\begin{defn} ${\cal X}\subseteq {\cal R}$ is ${\cal H}$-\textbf{Baire} if for every $[a,A]\neq \emptyset$, with $A\in {\cal H}$, there exists $\emptyset\neq[b,B]\subseteq[a,A]$, with $B\in{\cal H}$, such that $[b,B]\subseteq {\cal X}$ or $[b,B]\subseteq {\cal X}^c$. If for every $[a,A]\neq \emptyset$, with $A\in{\cal H}$, there exists $\emptyset\neq[b,B]\subseteq[a,A]$, with $B\in {\cal H}$, such that $[b,B]\subseteq {\cal X}^c$; we say that ${\cal X}$ is ${\cal H}$-\textbf{meager}.
\end{defn}

\noindent It is clear that if ${\cal X}\subseteq {\cal R}$ is ${\cal H}$-Ramsey then ${\cal X}$ is ${\cal H}$-Baire. Now, the notion corresponding to \emph{dense open} sets will be defined in this context: Given $A\in {\cal H}$ and ${\cal I}\subseteq{\cal AR}(A)$, we say that the sequence $({\cal D}_a)_{a\in{\cal I}}$, with ${\cal D}_a \subseteq {\cal H}$, $[a,C]\neq\emptyset$ for some $C\in{\cal D}_a$ and every $a\in {\cal I}$, has the \textbf{D--O property bellow} $A$ if for every $a\in{\cal I}$ the following hold:
\begin{enumerate}
\item $\forall B\in [a,A]\cap{\cal H}$ $\exists C\in {\cal D}_a$ $(C\leq B)$.
\item $[B\in {\cal D}_a\!\!\upharpoonright \!\!A$ $\wedge$ $(C\in[depth_B(a),B])\cap{\cal H}]$ $\Rightarrow$ $C\in {\cal D}_a$.
\end{enumerate}

\noindent The notion of selectivity is clear in this context:

\begin{defn} A coideal ${\cal H}\subseteq{\cal R}$ is \textbf{selective} if given $A\in {\cal H}$ and $(A_a)_{a\in{\cal I}}$, with ${\cal I}\subseteq{\cal AR}$, $A_a\in{\cal H}\!\!\upharpoonright \!\!A$ and $[a,A_a]\neq\emptyset$ for $a\in {\cal I}$, there exists $B\in{\cal H}\!\!\upharpoonright \!\!A$ such that $[a,B]\subseteq[a,A_a]$ for every $a\in{\cal I}\cap{\cal AR}(B)$.
\end{defn}

\noindent Now, it will be shown that selectivity implies the following property which will be useful in proving the main result of this work. The same name of the corresponding notion of coideals on $\mathbb{N}$ introduced by Farah will be used.

\begin{defn}
We say that ${\cal H}\subseteq {\cal R}$ is \textbf{semiselective} if given $A\in {\cal H}$, for every sequence $({\cal D}_a)_{a\in {\cal I}}$ with ${\cal I}\subseteq{\cal AR}(A)$, ${\cal D}_a\subseteq {\cal H}$ and with the D--O property below $A$, there exists $B\in {\cal H}\!\!\upharpoonright \!\!A$ such that $[depth_B(a),B]\cap{\cal H}\subseteq {\cal D}_a$ for every $a\in {\cal I}\cap {\cal AR}(B)$.
\end{defn}

\begin{prop} If ${\cal H}\subseteq{\cal R}$ is a selective coideal then ${\cal H}$ is semiselective.
\end{prop}

\noindent \textbf{Proof:} Given $A\in{\cal H}$, consider $({\cal D}_a)_{a\in {\cal I}}$ with ${\cal I}\subseteq{\cal AR}(A)$, ${\cal D}_a\subseteq {\cal H}$ and with the D--O property below $A$. For $a\in{\cal AR}(A)$, by (A.5) mod ${\cal H}$ there exists $B\in[a,A]\cap {\cal H}$ such that $[a,B]\cap{\cal H}\neq\emptyset$. By the D--O property, we can choose $A_a\in{\cal D}_a$ with $A_a\leq A$ and (again, by (A.5) mod ${\cal H}$) $[a,A_a]\neq\emptyset$. By selectivity, there exists $B\in {\cal H}\!\!\upharpoonright \!\!A$ such that $[a,B]\subseteq[a,A_a]$ for $a\in{\cal AR(B)}\cap {\cal I}$. But $[A_a\in {\cal D}_a\!\!\upharpoonright \!\!A$ $\wedge$ $(C\in[depth_B(a),B])\cap{\cal H}]$ $\Rightarrow$ $C\in {\cal D}_a$ (D--O property). Thus, $[depth_B(a),B])\cap{\cal H}]\subseteq{\cal D}_a$ for every $a\in{\cal AR(B)}\cap {\cal I}$. 

\qed

\vspace{.25cm}

\noindent The following is the version of theorem 1.6 from \cite{mij2} corresponding to this context and can be easily generalized to partitions in $n$ pieces:

\begin{thm}\label{ramsey2}
Suppose that ${\cal H}\subseteq {\cal R}$ is a selective coideal. Then, given a partition $f\colon {\cal AR}_2\to \{0,1\}$ and $A\in {\cal H}$, there exists $B\in {\cal H}\!\!\upharpoonright \!\!A$ such that $f$ is constant on ${\cal AR}_2(B)$.
\end{thm}

\noindent \textbf{Proof:} Let $f$ be the partition ${\cal AR}_2={\cal C}_0\cup {\cal C}_1$, and consider $A\in {\cal H}$. By (A.6)mod ${\cal H}$, for every $a\in {\cal AR}(A)$ we can define the nonempty
$${\cal D}_a=\{B\in [depth_A(a),A]\cap{\cal H}\colon f {\mbox { is constant on }}r_2([a,B])\}$$
if $a\in {\cal AR}_1(A)$, and ${\cal D}_a={\cal H}\!\!\upharpoonright \!\!A$ otherwise; which gives us a sequence with the D--O property below $A$. By selectivity (or the S--property), we have $\hat{B}\in {\cal H}\!\!\upharpoonright \!\!A$ such that $[depth_{\hat{B}}(a),\hat{B}]\subseteq{\cal D}_a$ for every $a\in {\cal AR}_1(\hat{B})$. Since $\hat{B}\in[depth_{\hat{B}}(a),\hat{B}]\subseteq{\cal D}_a$ for every $a\in {\cal AR}_1(\hat{B})$, there exists $i_a\in \{0,1\}$ such that $b\in {\cal C}_{i_a}$ if $b\in r_2([a,\hat{B}])$. Now, consider the partition $g\colon {\cal AR}_1\to \{0,1\}$ defined by $g(a)=i_a$ if $a\in {\cal AR}_1(\hat{B})$. By (A.6)mod ${\cal H}$, there exists $B\in [0,\hat{B}]\cap{\cal H}$ such that $g$ is constant on $r_1([0,B])={\cal AR}_1(B)$. But $B\leq \hat{B}\leq A$, so $B$ is as required. 

\qed

\vspace{0.4cm}

\noindent To give the local version of the corresponding Galvin lemma (or Nash-williams theorem) for selective coideals of ${\cal R}$, the following combinatorial forcing will be used: Fix ${\cal F}\subseteq {\cal AR}$. We say that $A\in {\cal H}$ \emph{\textbf{accepts}} $a\in {\cal AR}$ if $[a,A]=\emptyset$ or for every $B\in [a,A]\cap{\cal H}$ there exists $n\in \mathbb{N}$ such that $r_n(B)\in{\cal F}$. We say that $A$ \emph{\textbf{rejects}} $a$ if $[a,A]\neq\emptyset$ and no element of $[depth_A(a),A]\cap{\cal H}$ accepts $a$; and we say that $A$ \emph{\textbf{decides}} $a$ if $A$ either accepts or rejects $a$. This combinatorial forcing has the following properties:
\begin{enumerate}
\item If $A$ accepts $a$, then every $B\in {\cal H}\!\!\upharpoonright \!\!A$ accepts $a$.
\item If $A$ rejects $a$, then every $B\in {\cal H}\!\!\upharpoonright \!\!A$ rejects $a$, if $[a,B]\neq\emptyset$.
\item For every $A\in {\cal H}$ and every $a\in {\cal AR}(A)$ there exists $B\in [depth_A(a),A]\cap{\cal H}$ which decides $a$.
\item If $A$ accepts $a$ then $A$ accepts every $b\in r_{|a|+1}([a,A])$.
\item If $A$ rejects $a$ then there exists $B\in [depth_A(a),A]\cap{\cal H}$ which rejects every $b\in r_{|a|+1}([a,B])$.
\end{enumerate}

\noindent \textbf{Claim 1:} Given $A\in{\cal H}$, with $\mathcal{H}$ selective, there exists $B\in {\cal H} \!\!\upharpoonright \!\!A$ which decides every $b\in {\cal AR}(B)$.

\noindent \emph{Proof:} For every $a\in {\cal AR}(A)$ define
$${\cal D}_a=\{C\in[depth_A(a),A]\cap{\cal H}\colon C {\mbox{ decides }}a\}$$
Then $({\cal D}_a)_{a\in {\cal AR}(A)}$ has the D--O property, so there exists $B\in {\cal H}\!\!\upharpoonright \!\!A$ such that for every $a\in {\cal AR}(B)$ we have $[depth_B(a),B]\cap {\cal H}\subseteq {\cal D}_a$. Thus, $B$ decides every $a\in{\cal AR}(B)$. 

\qed

\begin{lem}\label{galvinlocal}
Given ${\cal F}\subseteq{\cal AR}$, a selective coideal ${\cal H}\subseteq \mathcal{R}$, and $A\in{\cal H}$, there exists $B\in {\cal H}\!\!\upharpoonright \!\!A$ such that one of the following holds:
\begin{enumerate}
\item ${\cal AR}(B)\cap{\cal F}=\emptyset$, or
\item $\forall C\in {\cal H}\!\!\upharpoonright \!\!B$ $(\exists \ n\in \mathbb{N})$ $(r_n(C)\in{\cal F})$.
\end{enumerate}
\end{lem}

\noindent \textbf{Proof:} consider $B$ as in claim 1. If $B$ accepts $\emptyset$ part (2) holds. Assume that $B$ rejects $\emptyset$ and for $a\in {\cal AR}(B)$ define
$${\cal D}_a=\{C\in[depth_B(a),B]\cap{\cal H}\colon C {\mbox{ rejects every }}b\in r_{|a|+1}([a,C])\}$$if $B$ rejects $a$, and ${\cal D}_a={\cal H}\!\!\upharpoonright \!\!B$ otherwise. So, $({\cal D}_a)_{a\in {\cal AR}(A)}$ has the D--O property bellow $A$. Then we have $\hat{B}\in {\cal H}\!\!\upharpoonright \!\!B$ such that $[depth_{\hat{B}}(a),\hat{B}]\cap {\cal H}\subseteq {\cal D}_a$ for every $a\in {\cal AR}(\hat{B})$. By induction on the lenght, $\hat{B}$ rejects every $a\in {\cal AR}(\hat{B})$, hence no element of ${\cal AR}(\hat{B})$ is in ${\cal F}$. Thus, part (1) holds. 

\qed

\begin{thm}\label{baire-ramsey}
If $\mathcal{H}\subseteq\mathcal{R}$ is a selective coideal then ${\cal X}\subseteq{\cal R}$ is ${\cal H}$--Ramsey iff ${\cal X}$ is ${\cal H}$--Baire
\end{thm}

\noindent \textbf{Proof:} Let ${\cal X}$ be a ${\cal H}$--Baire subset of ${\cal R}$ and consider $A\in {\cal H}$. As before, we only proof the result for $[\emptyset,A]$ without loss of generality. For $a\in{\cal AR}(A)$ define
$${\cal D}_a=\{B\in[depth_A(a),A]\cap{\cal H}\colon [a,B]\subseteq {\cal X} {\mbox { or }} [a,B]\subseteq {\cal X}^c$$ $${\mbox { or }} ([a,C]\not\subseteq{\cal X} {\mbox { and }} [a,C]\not\subseteq{\cal X}^c \ \forall C\in[a,B] )\}$$
Then $({\cal D}_a)_a$ has the D--O property bellow $A$. Let $\hat{B}\in{\cal H}\!\!\upharpoonright \!\!A$ such that, for $a\in{\cal AR}(\hat{B})$, $[depth_{\hat{B}}(a),\hat{B}]\cap{\cal H}\subseteq {\cal D}_a$ . Let ${\cal F}_0=\{a\in{\cal AR}(A)\colon [a,\hat{B}]\subseteq {\cal X}\}$ and ${\cal F}_1=\{a\in{\cal AR}(A)\colon [a,\hat{B}]\subseteq {\cal X}^c\}$. By applying lemma \ref{galvinlocal} to ${\cal F}_0$ (or to ${\cal F}_1$) and $\hat{B}$, we obtain $B\in{\cal H}\!\!\upharpoonright \!\!\hat{B}$ such that $[\emptyset,B]\subseteq{\cal X}$ (or $[\emptyset,B]\subseteq{\cal X}^c$) or ${\cal AR}(B)\cap({\cal F}_0\cup{\cal F}_1)=\emptyset$. The latter case is not possible: since ${\cal X}$ is ${\cal H}$--Baire, there exists $\emptyset\neq[b,C]\subseteq [\emptyset,B]$ such that $[b,C]\subseteq{\cal X}$ or $[b,C]\subseteq{\cal X}^c$. By (A.5) mod ${\cal H}$, we can suppose that $C\in[b,C]\subseteq[b,\hat{B}]$, and since $\hat{B}\in{\cal D}_b$, we conclude that $b\in{\cal F}_0\cup{\cal F}_1$. The reverse implication is obvious. 

\qed

\vspace{0.4cm}

\noindent Now, we give one more application of lemma \ref{galvinlocal}. Recall that the \emph{metric open} subsets of ${\cal R}$ are of the form
$$[b]=\{A\in{\cal R}\colon b\sqsubseteq A\}$$where $b\sqsubseteq A$ means $\exists n\in\mathbb{N}(r_n(A)=b)$.

\begin{thm}\label{abiertos}
Suppose that $\mathcal{H}\subseteq\mathcal{R}$ is a selective coideal. Then the metric open subsets of ${\cal R}$ are ${\cal H}$-Ramsey.
\end{thm}
\noindent \textbf{Proof:} Let ${\cal X}\subseteq{\cal R}$ metric open and consider $[\emptyset,A]$ with $A\in {\cal H}$. Define, for every $a\in{\cal AR}$
$${\cal D}_a=\{B\colon [a,B]\subseteq {\cal X} {\mbox { or }}[a,B]\subseteq{\cal X}^c{\mbox { or }}\forall C\leq B([a,C]\not\subseteq {\cal X} {\mbox { and }}[a,C]\not\subseteq {\cal X}^c)\}$$
Thus, $\{{\cal D}_a\}_a$ has the D-O property below $A$. So, let $B\in{\cal H}\!\!\upharpoonright \!\!A$ be such that $[depth_A(a),B]\cap{\cal H}\subseteq {\cal D}_a$ for every $a\in{\cal AR}(B)$. Now, define $h\colon {\cal AR}\to \{0,1,2\}$ by
$$h(a)=\left\{
\begin{array}{cc}
0 & {\mbox { if }}[a,B]\subseteq{\cal X}\\
1 & {\mbox { if }}[a,B]\subseteq{\cal X}^c\\
2 & {\mbox { otherwise}}
\end{array}
\right.$$

\noindent If $h(a)=2$, by restricting $h$ to ${\cal AR}_{|a|+1}$ we obtain (by (A.6) mod${\cal H}$) $C_a\in[depth_B(a),B]\cap{\cal H}$ such that $h$ is constant on $r_{|a|+1}([a,C_a])$. Furthermore, that constant is $2$ since $h(a)=2$ and $B\in{\cal D}_a$.

\begin{clm}
If $h(\emptyset)=2$ then there exists $C\in {\cal H}\!\!\upharpoonright \!\!B$ such that $\forall b\in{\cal AR}(C)(h(c)=2)$.
\end{clm}
\noindent \textbf{Proof:}(of the claim) Define for $b\in{\cal AR}$, $C_b$ as before if $[b,B]\neq\emptyset$ and $C_b=B$ otherwise. Then $\{C_b\}_b$ has the D-O property. Let $C\in{\cal H}\!\!\upharpoonright \!\!B$ be such that $[depth_C(b),C]\cap{\cal H}\subseteq {\cal C}_a$ for every $b\in{\cal AR}(C)$. Suppose that $h(b)\neq2$ for some $b=r_{|b|}(\hat{B})\in{\cal AR}(B)$ and choose it with minimal depth in $C$. Thus, $b\neq\emptyset$ since $h(\emptyset)=2$. Let $b'=r_{|b|-1}(\hat{B})$. Then $h(b')=2$, but $b\in r_{|b'|+1}([b',C])\subseteq r_{|b'|+1}([b',C_{b'}])$ and hence $h(b)=2$ (see the paragraph before the claim). This is a contradiction, and the claim is proved. 

\qed

\vspace{0.4cm}

\noindent Now it will be shown that $h(\emptyset)<2$. Suppose that $h(\emptyset)=2$ and $C$ is as in the claim. Then $[\emptyset,C]\not\subseteq{\cal X}$ and $[\emptyset,C]\not\subseteq{\cal X}^c$. Consider $\hat{C}\in {\cal X}\!\!\upharpoonright \!\!C$. Since ${\cal X}$ is metric open, there exists $b\in {\cal AR}$ such that $b\sqsubseteq B$ and $[b]\subseteq {\cal X}$, i. e., $h(b)=0$, which is a contradiction (by the claim). This completes the proof of theorem. 

\qed

\section{The Souslin operation}

\noindent The goal of this section is to show that the family of ${\cal H}$--Ramsey subsets of ${\cal R}$ is closed under the Souslin operation when ${\cal H}$ is a selective coideal.

\begin{lem} If ${\cal H}\subseteq{\cal R}$ is a selective coideal of $\mathcal{R}$ then the families of ${\cal H}$--Ramsey and ${\cal H}$--Ramsey null subsets of ${\cal R}$ are closed under countable union.
\end{lem}

\noindent \textbf{Proof:} Fix $A\in{\cal H}$. The proof will be made for $[\emptyset,A]$ without loss of generality. Suppose that $({\cal X}_n)_{n\in\mathbb{N}}$ is a sequence of ${\cal H}$--Ramsey null subsets of ${\cal R}$. Define for $a\in{\cal AR}(A)$
$${\cal D}_a=\{B\in[a,A]\cap{\cal H}\colon [a,B]\subseteq {\cal X}_n^c \ \ \forall n \leq |a| \}$$
Then $({\cal D}_a)_a$ has the D--O property bellow $A$, so let $B\in{\cal H}\!\!\upharpoonright \!\!A$ be such that $[depth_B(a),B]\cap{\cal H}\subseteq {\cal D}_a$ for all $a\in {\cal AR}(B)$. Thus, $[\emptyset, B]\subseteq \bigcap {\cal X}_n^c$ (since $B\in[depth_B(a),B]\cap{\cal H}$ for every $a\in{\cal AR}(B)$). Now, suppose that $({\cal X}_n)_{n\in\mathbb{N}}$ is a sequence of ${\cal H}$--Ramsey subsets of ${\cal R}$ and consider $\emptyset\neq[a,A]$. If there exists $B\in{\cal H}\!\!\upharpoonright \!\!A$ such that $[a,B]\subseteq {\cal X}_n$ for some $n$, we are done. Otherwise, using an argument similar to the one above, we prove that $\bigcup {\cal X}_n$ is ${\cal H}$--Ramsey null. 

\qed

\vspace{0.4cm}

\noindent Recall that given a set $X$, two subsets $A,B$ of $X$ are "\emph{compatibles}" with respect to a family ${\cal F}$ of subsets $X$ if there exists $C\in {\cal F}$ such that $C\subseteq A\cap B$. And ${\cal F}$ is \emph{M-like} if for ${\cal G}\subseteq {\cal F}$ with $|{\cal G}|<|{\cal F}|$, every member of ${\cal F}$ which is not compatible with any member of ${\cal G}$ is compatible with $X\setminus \bigcup {\cal G}$. A $\sigma$-algebra ${\cal A}$ of subsets of $X$ together with a $\sigma$-ideal ${\cal A}_0\subseteq{\cal A}$ is a \emph{Marczewski pair} if for every $A\subseteq X$ there exists $\Phi(A)\in {\cal A}$ such that $A\subseteq \Phi(A)$ and for every $B\subseteq \Phi(A)\setminus A$, $B\in{\cal A}\rightarrow B\in {\cal A}_0$. The following is a well known fact:

\begin{thm}[Marczewski]\label{marcz}
Every $\sigma$-algebra of sets which together with a $\sigma$-ideal is a Marczeswki pair, is closed under the Souslin operation. 
\end{thm}
\qed

\noindent Denote $Exp({\cal H})=\{[n,A]\colon n\in \mathbb{N}$, $A\in {\cal H}\}$. ${\cal H}$ selective coideal of $\mathcal{R}$.

\begin{prop}\label{mlike}
If $|{\cal H}|=2^{\aleph_0}$, then the family $Exp({\cal H})$ is $M$-like.
\end{prop}

\noindent \textbf{Proof:} Consider ${\cal B}\subseteq Exp({\cal H})$ with $|{\cal B}|<|Exp({\cal H})|=2^{\aleph_0}$ and suppose that $[a,A]$ is not compatible with any member of ${\cal B}$, i. e. for every $B\in {\cal B}$, $B\cap [a,A]$ does not contain any member of $Exp({\cal H})$. We claim that $[a,A]$ is compatible with ${\cal R}\smallsetminus \bigcup {\cal B}$. In fact:

\noindent Since $|{\cal B}|<2^{\aleph_0}$, $\bigcup {\cal B}$ is ${\cal H}$-Baire  (it is ${\cal H}$-Ramsey). So, there exist $[b,B]\subseteq [a,A]$ with $B\in {\cal H}$ such that:
\begin{enumerate}
\item $[b,B]\subseteq \bigcup {\cal B}$ or
\item $[b,B]\subseteq {\cal R}\smallsetminus \bigcup {\cal B}$
\end{enumerate}

\noindent(1) is not possible because $[a,A]$ is not compatible with any member of ${\cal B}$. And (2) says that $[a,A]$ is compatible with ${\cal R}\smallsetminus \bigcup {\cal B}$ 

\qed

\vspace{0.4cm}

\noindent Proposition \ref{mlike} says that the family of ${\cal H}$--Ramsey subsets of ${\cal R}$ together with the family of ${\cal H}$--Ramsey null subsets of ${\cal R}$ is a Marczewski pair (see section 2 of \cite{pawl}). Thus, by theorem \ref{marcz}, we have the result:

\begin{thm}\label{souslin}
The family of ${\cal H}$--Ramsey subsets of ${\cal R}$ is closed under the Souslin operation. 
\end{thm}

\qed

\noindent As a consequence of theorems \ref{abiertos} and \ref{souslin}, the following holds:

\begin{thm}\label{analiticos}
Suppose that $\mathcal{H}\subseteq\mathcal{R}$ is a selective coideal. Then the analitic subsets of ${\cal R}$ are ${\cal H}$-Ramsey. 
\end{thm}

\qed

\section{Examples}\label{examples}
The goal of this section is to give examples of topological Ramsey spaces for which the previous results are ilustrated. An example from the Ellentuck's topological Ramsey space $(\mathbb{N}^{[\infty]},\subseteq,r_n)$ as defined in section \ref{topo} is the following: fix $(x_n)_n\subseteq\mathbb{N}^{\mathbb{N}}$ and $x$ a cluster point of $(x_n)_n$. Define 
$$\mathcal{H}=\{A\in\mathbb{N}^{[\infty]}\colon x {\mbox{ is cluster point of }(x_n)_{n\in A}}\}$$ 
Then, $\mathcal{H}$ is a coideal (in our context and in the sense of the known notion of coideal). If $x_n$ is borel for every $n\in \mathbb{N}$, then $\mathcal{H}$ is selective and $\mathbb{N}^{[\infty]}\setminus \mathcal{H}$ is analitic, that is to say, $\mathcal{H}$ is $\Pi_1^1$.

\vspace{.25cm}

\noindent Another example: Fix $k\in \mathbb{N}$. Given $p\colon \mathbb{N}\to \{0,1,\dots,k\}$, denote $supp(p)=\{n\colon p(n)\neq 0\}$ and $rank(p)$ the image set of $p$.
Consider the set 
$$FIN_k=\{p\colon\mathbb{N}\to \{0,1,\dots\}\colon |supp(p)|<\infty {\mbox { and }} k\in rank(p)\}$$  
we say that $X=(x_n)_{n\in \mathcal{I}}\subseteq FIN_k$, with $\mathcal{I}\in \mathcal{P}(\mathbb{N})$ is a \emph{basic block sequence} if 
\begin{center}
$n<m\Rightarrow $max($supp(x_n))<$min($supp(x_m))$ 
\end{center}
For infinite sequences we assume that $\mathcal{I}=\mathbb{N}$. Define $T\colon FIN_k \to FIN_{k-1}$ by 
\begin{center}
$T(p)(n)=$max$\{p(n)-1,0\}$
\end{center}
For $j\in\mathbb{N}$, $T^{(j)}$ is the $j$-th iteration of $T$. Given a basic block sequence $X=(x_n)_{n\in\mathcal{I}}$ we define $[X]\subseteq FIN_k$ as the set which elements are of the form
$$T^{(j_0)}(x_{n_0})+T^{(j_1)}(x_{n_1})+\cdots+T^{(j_r)}(x_{n_r})$$
with $n_0<n_1<\cdots<n_r\in \mathcal{I}$, $j_0<j_1<\cdots<j_r\in\{0,1,\dots,k\}$, and $j_i=0$ for some $i\in\{0,1,\dots,r\}$. Denote $FIN_k^{[\infty]}$, the set of infinite basic block sequences, for $A$, $B\in FIN_k^{[\infty]}$, define
$$A\leq B \Leftrightarrow A\subseteq [B]$$and $r_n(A)=$ "the first $n$ elements of $A$". Then $(FIN_k^{[\infty]}, \leq, r)$ is a topological Ramsey space (see \cite{todo}). Furthermore, we have the following well known result:

\begin{thm}[Gowers]\label{gowers}
Given an integer $n>0$ and $f\colon FIN_k\to \{0,1,\dots,r-1\}$, there exists $A\in FIN_k^{[\infty]}$ such that $f$ is constant on $[A]$.  
\end{thm}
\qed

\noindent For $k=1$, the previous theorem reduces to the famous Hindman's theorem (\cite{hindman}). Assuming \textbf{CH}, we define a well order $(\mathcal{P}(FIN_k),<)$, an for a fixed $X\subseteq FIN_k$ we find $A_X\in FIN^{[\infty]}$ such that
\begin{enumerate}
\item $\mathcal{AR}(A_X)\subseteq X$ or $\mathcal{AR}(A_X)\subseteq X^c$.
\item $X<Y\Rightarrow A_X\leq^*A_Y$
\end{enumerate}
Where "$A\leq^*B$" means that $A\leq B$ "from some $n$ on". Suppose that we have defined $A_Y$ for every $Y<X$. We only have to consider the  case in which $X$ is limit. If we have already $A_{Y_0}\geq A_{Y_1}\geq \cdots$ for the predecesors $Y_0,Y_1,\dots$ of $X$, we can choose $a_0\in A_{Y_0}$, $a_1\in A_{Y_1}$, $\dots$ such that $a_0<a_1<\cdots$. Then, $A=a_0^{\smallfrown}a_1^{\smallfrown}\cdots$ satisfies (2). Now, if 
$$\mathcal{AR}(A)=(X\cap \mathcal{AR}(A))\cup ((X^c\cap \mathcal{AR}(A)))$$we can find, by theorem \ref{gowers}, $A_X\leq A$ which satisfies (1). It is clear that $A_X$ satisfies (2) too. This completes the construction.

Now, the coideal:

$$\mathcal{H}=\{B\in FIN_k^{[\infty]}\colon \exists X\subseteq FIN_k (A_X\leq B)\}$$   
It is clear that $\mathcal{H}$ satisfies (1) and (2) from the definition of coideal. Now, by theorem \ref{gowers}, and the previous construction, given $B\in\mathcal{H}$ and $f\colon [B]\to \{0,1\}$ there exists $B'\in \mathcal{H}$ such that $f$ is constant on $[B']$. This is, $(A.6)$mod $\mathcal{H}$ holds. The selectivity of $\mathcal{H}$ is also a consequence of the construction of the $A_X$'s (which is strongly based on \textbf{CH}, of course. See \cite{blass}). The previous construction can be done on any topological Ramsey space in a similar way, under the assumption of \textbf{CH} or the Martin's axiom. That is to say, What is given above is a scheme of examples.

\end{document}